\documentclass{gtart}

\usepackage{amsfonts}
\usepackage{amsmath}
\usepackage{amssymb}
\usepackage{amstext}
\usepackage{amsthm}       
\usepackage{xspace}
\usepackage{graphicx}
\usepackage{psfrag}
\usepackage{epic,eepic}

\newcommand{\R}{\mathbb R}

\newcommand{\Q}{\mathbb Q}
\newcommand{\Z}{\mathbb Z} 

\newtheorem{theorem}{Theorem}[section]
\newtheorem{lemma}[theorem]{Lemma}
\newtheorem{cor}[theorem]{Corollary}
\newtheorem{prop}[theorem]{Proposition}

\theoremstyle{definition}
\newtheorem{dfn}[theorem]{Definition}
\newtheorem*{rem}{Remark}

\begin{document}

\title{On tight contact structures with negative maximal twisting 
number on small Seifert manifolds}
\author{Paolo Ghiggini}
\begin{abstract}
We study some properties of transverse contact structures 
on small Seifert manifolds, and we apply them to the classification 
of tight contact structures on a family of small Seifert manifolds.
\end{abstract}
\maketitle

\section{Introduction}
\begin{figure}
\begin{center}
\psfrag{0}{$e_0$}
\psfrag{a}{$-\frac{1}{r_1}$}
\psfrag{b}{$-\frac{1}{r_2}$}
\psfrag{c}{$-\frac{1}{r_3}$}
\includegraphics[width=10cm]{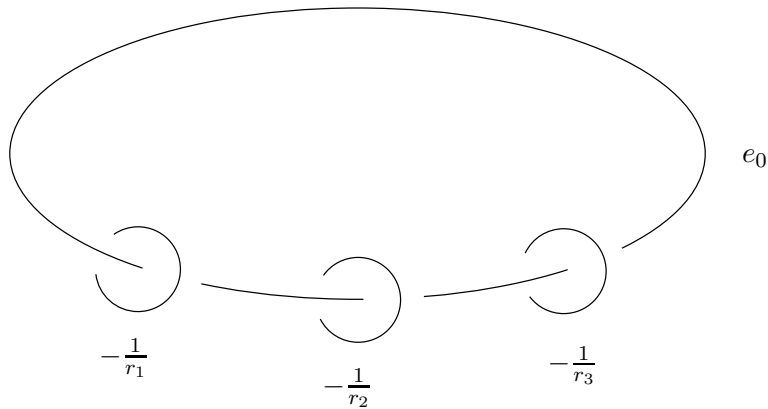}
\end{center}
\caption{Rational surgery diagram for the Seifert fibred 3--manifold 
$M(e_0;r_1,r_2,r_3)$.}
\label{1.fig}
\end{figure}

In this article $M(e_0; r_1,r_2,r_3)$ --- with $e_0 \in \Z$ and $r_i \in (0,1) \cap 
\Q$ --- will denote the 
$3$--manifold specified by the surgery diagram in Figure \ref{1.fig}.
It is well known that $M(e_0; r_1,r_2,r_3)$ carries a Seifert fibration over
$S^2$ with three singular fibres corresponding to the three small 
unknots in the surgery diagram. The manifolds belonging to this family
 will be called {\em small Seifert manifolds}. 

The classification of tight contact structures on small Seifert 
manifolds has been the object of intense study in the last few years.
The generic case, when $e_0 \neq -1, -2$, was settled in \cite{hwu} and 
\cite{gls:1}, and a large family of manifolds with $e_0=-1$ was studied
in \cite{gls:2}. The goal of this article is the classification of 
tight contact structures on some small Seifert manifolds with $e_0=-2$.
Such results are useful in symplectic cut-and-paste operations like
the generalised symplectic rational blow-down \cite{gay-stipsicz:1}.

The main invariant in the classification of tight contact structures 
on Seifert manifolds is the maximal twisting number. Let $L$ be a 
regular fibre for the Seifert fibration on $M= M(e_0; r_1,r_2,r_3)$, and 
let ${\mathcal S}$ be the set of isotopies $\varphi \colon [0,1] \times M \to M$ such 
that $\varphi_0$ is the identity and $\varphi_1(L)$ is a Legendrian curve. $L$ has a 
distinct framing induced by the Seifert fibration, so we can transport 
this framing to $\varphi_1(L)$. We denote by $L_{\varphi}$ the framed curve $\varphi_1(L)$ 
with the framing induced by $\varphi$. As a Legendrian curve, $\varphi_1(L)$ has also
a framing induced by the contact structure. We define the {\em twisting 
number} $tb(L_{\varphi})$ as the difference between the contact framing and 
the framing induced by $\varphi$.

\begin{dfn}
For any contact structure $\xi$ on $M$, we define 
the {\em maximal twisting number} of $\xi$ as
\[
t(\xi)= \max_{\varphi \in  {\mathcal S}} \min \{ tb(L_{\varphi}), 0 \}.
\]
%where ${\mathcal S}$ is the set of all Legendrian curves $L \subset M$ 
%which are smoothly isotopic to a regular fibre.
\end{dfn}

We can see ${\mathcal S}$ as the universal cover of the space of 
Legendrian curves isotopic to a regular fibre ({\em vertical 
Legendrian curves} from now on). However we  
would prefer to see the twisting number as a function on the space of
vertical Legendrian curves, not on its universal cover. This is the 
case when the framing on $L_{\varphi}$ is independent of $\varphi$, and happens in 
the manifolds we are interested in. In fact if two isotopies induce 
different framings on the same vertical Legendrian curve, then the 
twisting number can be made arbitrarily big, so the contact 
structures has $t=0$. Moreover, if this happens for one contact 
structure, then it happens for all.

The first result of this article is a necessary and sufficient 
condition for the existence of tight contact structures with negative 
maximal twisting number on small Seifert manifolds with $e_0=-1$.  
Necessary and sufficient conditions 
for small Seifert manifolds with $e_0 \neq -1$ were given in \cite{hwu}: 
tight contact structures with negative maximal twisting number exist
on a small Seifert manifold with $e_0 \neq -1$ if and only if $e_0$ is 
negative. Our condition is the following:

\begin{theorem} \label{principale}
Let $M(e_0; r_1, r_2, r_3)$ be a small Seifert manifold with
$e_0=-1$. Then the following facts are equivalent:
\begin{enumerate}
\item $M(e_0; r_1, r_2, r_3)$ carries a tight contact structure $\xi$ with 
$t(\xi)<0$
\item there exist integer numbers $p_1$, $p_2$, $p_3$, and $q$ with $q>0$
such that 
\begin{itemize}
\item $(p_i, q)=1$ and $\frac{p_i}{q} < -r_i$,
\item if $\frac{p'}{q'} \in (\frac{p_i}{q}, -r_i)$, then $q' > q$
\end{itemize}
\item $M(e_0, r_1,r_2,r_3)$ carries a contact structure transverse to the 
Seifert fibration.
\end{enumerate}
\end{theorem}

For each of the three rational numbers $r_1$, $r_2$, $r_3$,
we can write
\begin{equation*}
- \frac{1}{r_i}= [a^{(i)}_0,a^{(i)}_1,\ldots,a^{(i)}_{m_i}] =
a_0^{(i)} - \cfrac{1}{a_1^{(i)} - \cfrac{1}{\ddots -
  \cfrac{1}{a_{m_i}^{(i)}}}},\quad i=1,2,3,
\end{equation*}
for some uniquely determined integer coefficients 
\[ 
a^{(i)}_0,\cdots,a^{(i)}_{m_i}\leq -2, \quad i=1,2,3.
\]
We define
\[ T(r_i)=  \prod_{k=0}^{m_i} (a_k^{(i)}+1). \]

\begin{figure}

\setlength{\unitlength}{1pt}

\begin{picture}(420,180)(-210,-250)
% bottom

% bottom center

\put(-10,-85){\tiny{$-2$}}

\put(0,-110){\arc{40}{9.04}{11.95}}

\put(0,-110){\arc{40}{5.9}{7.25}}

\put(0,-110){\arc{40}{7.45}{8.8}}

% bottom left

\put(-40,-85){\tiny{$a^{(1)}_1$}}

\put(-35,-110){\arc{40}{9.04}{11.95}}

\put(-35,-110){\arc{40}{5.9}{8.8}}

\put(-65,-85){\tiny{$a^{(1)}_2$}}

\put(-70,-110){\arc{40}{11}{11.95}}

\put(-70,-110){\arc{40}{5.9}{7.85}}

\put(-83,-110){$\cdots$}

\put(-105,-85){\tiny{$a^{(1)}_{m_1-1}$}}

\put(-85,-110){\arc{40}{9.04}{11}}

\put(-85,-110){\arc{40}{7.85}{8.8}}

\put(-130,-85){\tiny{$a^{(1)}_{m_1}$}}

\put(-120,-110){\arc{40}{5.9}{11.95}}

% bottom right

\put(30,-85){\tiny{$a^{(3)}_1$}}

\put(35,-110){\arc{40}{9.04}{11.95}}

\put(35,-110){\arc{40}{5.9}{8.8}}

\put(55,-85){\tiny{$a^{(3)}_2$}}

\put(70,-110){\arc{40}{9.04}{11}}

\put(70,-110){\arc{40}{7.85}{8.8}}

\put(72,-110){$\cdots$}

\put(85,-85){\tiny{$a^{(3)}_{m_3-1}$}}

\put(85,-110){\arc{40}{11}{11.95}}

\put(85,-110){\arc{40}{5.9}{7.85}}

\put(115,-85){\tiny{$a^{(3)}_{m_3}$}}

\put(120,-110){\arc{40}{9.04}{15.1}}

% bottom lower

\put(23,-145){\tiny{$a^{(2)}_1$}}

\put(0,-145){\arc{40}{4.33}{7.25}}

\put(0,-145){\arc{40}{7.45}{10.38}}

\put(23,-175){\tiny{$a^{(2)}_2$}}

\put(0,-180){\arc{40}{4.33}{6.28}}

\put(0,-180){\arc{40}{9.42}{10.38}}

\put(-2,-190){$\vdots$}

\put(23,-205){\tiny{$a^{(2)}_{m_2-1}$}}

\put(0,-195){\arc{40}{7.45}{9.42}}

\put(0,-195){\arc{40}{6.28}{7.25}}

\put(23,-230){\tiny{$a^{(2)}_{m_2}$}}

\put(0,-230){\arc{40}{4.33}{10.38}}

\end{picture}
\caption{Integer surgery presentation of $M(e_0; r_1, r_2, r_3)$.}
\label{link.fig}
\end{figure}
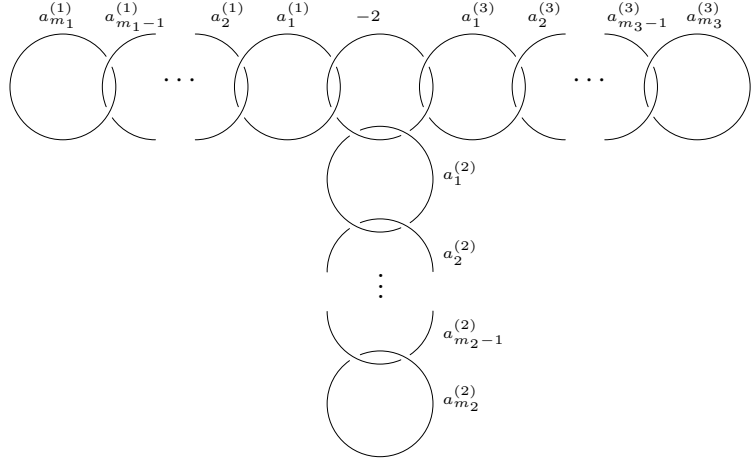

With the help of Theorem \ref{principale} we can classify tight 
contact structures on small Seifert manifolds with $e_0=-2$ which are 
$L$--spaces. An $L$--spaces was originally defined as a rational 
homology sphere $Y$ for which ${\rm rk} \widehat{HF}(Y) = |H_1(Y, \Z)|$.
However, thank to \cite[Theorem 1.1]{lisca-stipsicz:6}, small Seifert 
manifolds $M$ with $e_0=-2$ which are $L$--spaces can be characterised 
as those for which $-M$ carries no contact structures transverse to 
the Seifert fibration.

\begin{theorem} \label{classificazione}
Let $M(e_0; r_1, r_2, r_3)$ with $e_0 = -2$ be an $L$--space. Then all 
tight contact structures on $M(e_0; r_1, r_2, r_3)$ have maximal 
twisting number $-1$ and are Stein fillable. $M(e_0; r_1, r_2, r_3)$ 
admits exactly $T(r_1)T(r_2)T(r_3)$ isotopy classes of tight 
contact structure, which are distinct by their 
${\rm Spin}^c$--structures, and are filled by the Stein manifolds 
described by Legendrian surgery on all possible Legendrian 
realisations of the link in Figure \ref{link.fig}. 
\end{theorem} 

Theorem \ref{classificazione} confirms the conjecture that, 
on small Seifert manifolds which are $L$--spaces, all tight contact 
structure have non trivial Ozsv\'ath--Szab\'o and are distinguished by
the induced ${\rm Spin}^c$--structure \cite[Conjecture 1.2]{gls:2}.

\subsection*{acknowledgements}
During the final preparation of this article the author was supported 
by CIRGET and by the Chaire de Recherche 
du Canada en alg\`ebre, combinatoire et informatique math\'ematique 
de l'UQAM.

We warmly thank Paolo Lisca for sharing a preliminary version of
his work on transverse contact structures on Seifert manifolds
\cite{lisca-matic:transverse}, Andr\'as Stipsicz for his 
encouragement and for suggesting the statement of Theorem 
\ref{classificazione}, and Hao Wu for sharing Figure \ref{link.fig}.

Some of the results in this paper have been independently obtained also
by Patrick Massot. We thank him for carefully reading the first 
version of this article and for suggesting several improvements.

\section{Decomposition of negative twisting contact manifolds}
In this section we prove constraints on the maximal twisting number, 
in particular proving a necessary condition for the existence of
tight contact structures with negative maximal twisting number.
The reader is assumed to be familiar with convex surfaces theory
\cite{giroux:1} and bypasses \cite{honda:1}.

Let $V_i$ be a tubular neighbourhood of the singular fibre $F_i$ for
$i=1,2,3$. 
We identify $- \partial (M \setminus V_i)$ with $\R^2 / \Z^2$ so that 
$\binom{0}{1}$ is the direction of the regular fibres, 
and the meridian of $V_i$ has slope $- r_i$ in $- \partial (M \setminus V_i)$.

$M \setminus (V_1 \cup V_2 \cup V_3)$ is diffeomorphic to $\Sigma \times S^1$ where $\Sigma$ is a
pair of pants, and sometime it will be useful to consider also a 
second set of coordinates 
on $- \partial (M \setminus V_i)$ coming from the product structure such that 
$\binom{1}{0}$ is the direction of the section $\Sigma \times \{ 1 \}$. We choose
the diffeomorphism between $M \setminus (V_1 \cup V_2 \cup V_3)$ and $\Sigma \times S^1$ so that 
boundary slopes $s_1$, $s_2$ and $s_3$ in the old bases will correspond to 
boundary slopes $s_1'=s_1$, $s_2'=s_2$, and $s_3'=s_3-e_0$ in the bases coming 
from the product structure.
 
\begin{prop} \label{decomposizione}
Let $M=M(e_0; r_1, r_2, r_3)$ be a small Seifert manifold with integer 
Euler class $e_0=-1$ or $e_0=-2$, and let $\xi$ be a tight contact 
structure on $M$ with maximal twisting number $t(\xi)=-q<0$. Then for 
$i=1,2,3$ there exist tubular neighbourhoods $U_i$ of the singular 
fibres $F_i$ such that $M \setminus (U_1 \cup U_2 \cup U_3)$ has minimal convex boundary 
with  slopes $s(- \partial (M \setminus U_i))= \frac{p_i}{q}$ with $q>0$  satisfying
\begin{enumerate}
\item $(p_i, q)=1$ and $\frac{p_i}{q} < -r_i$, 
\item if $\frac{p'}{q'} \in (\frac{p_i}{q}, -r_i)$, then $q' < q$.
\end{enumerate}
\end{prop}
\begin{proof}
Let $L$ be a vertical Legendrian curve with twisting number 
$tb(L)=-q$. Isotope the Seifert fibration so that $L$ becomes a 
regular fibre and the singular fibres $F_i$ become Legendrian, then 
take standard neighbourhoods $V_i$ of $F_i$. We can make the twisting 
numbers of $F_i$ negative and as big as we wish in absolute value, 
therefore making the slopes of $- \partial (M \setminus V_i)$ arbitrarily close to $-r_i$.
 By repeatedly attaching the bypasses coming from convex vertical 
annuli with Legendrian boundary between $L$ and $- \partial (M \setminus V_i)$ 
(Imbalance Principle \cite[Proposition 3.17]{honda:1}),
we obtain tubular neighbourhoods $U_i$ of $F_i$ containing $V_i$ such that 
$- \partial (M \setminus U_i)$ has slope $\frac{p_i}{q}$. The numbers $\frac{p_i}{q}$ are 
uniquely determined by being the first ones which have denominator 
not greater 
than $q$ in the shortest path in the Farey Tessellation from the slope
of $- \partial (M \setminus V_i)$ to infinity. In particular, if $t(\xi)=-1$, then $q=1$,
$\frac{p_i}{q}=[-r_i]$, and the properties of $\frac{p_i}{q}$ follow 
immediately.

Assume now $t(\xi)< -1$, so that $q>1$. The property $(p_i, q)=1$ follows 
from $t(\xi)=-q$ because, if the fraction 
$\frac{p_i}{q}$ could be reduced for some $i$, then the twisting number
of a vertical 
Legendrian ruling curve of $- \partial (M \setminus U_i)$ would be greater that $-q$. 
Since the vertical Legendrian ruling curves 
of $- \partial (M \setminus U_i)$ are smoothly isotopic to regular fibres, this would 
contradict $t(\xi)=-q$. Since $- \partial (M \setminus U_i)$ is obtained from 
$- \partial (M \setminus V_i)$ by attaching vertical bypasses, and the attachment of a 
vertical bypasses decreases the slope, we have $\frac{p_i}{q} < -r_i$ for 
$i=1,2$.

We prove point $2$ by contradiction. Assume there is a rational 
number $\frac{p'}{q'} \in (\frac{p_i}{q}, -r_i)$ with $q' \leq  q$ for some 
$i$. If $q'=q$, then 
$\frac{p_i+1}{q} \in (\frac{p_i}{q}, -r_i)$. By the following 
Algebraic Lemma \ref{algebrico}, there is a fraction 
$\frac{p'}{q-1} \in (\frac{p_i}{q}, \frac{p_i+1}{q}]$, 
therefore we can assume $q' <  q$. 

By \cite[Lemma 2.15]{etnyre-honda:3} there is a neighbourhood $U_i'$ of 
the singular fibre $F_i$ such that $- \partial (M \setminus U_i')$ has slope 
$\frac{p'}{q'}$. If we put $- \partial (M \setminus U_i')$ in standard form, a vertical 
Legendrian ruling curve will be a Legendrian curve with twisting 
number $-q' > -q$ smoothly isotopic to a regular fibre. This 
contradicts the hypothesis $t(\xi)=-q$. 
\end{proof}

\begin{lemma} \label{algebrico}
For any rational number represented by a fraction $\frac pq$ with 
$q>1$ there exists an integer number $p'$ such that 
$\frac{p'}{q-1} \in [\frac{p}{q}, \frac{p+1}{q}]$. Moreover, if
$\frac{p}{q}$ is a reduced fraction, then $\frac{p'}{q-1} \neq 
\frac{p}{q}$.
\end{lemma}
\begin{proof}
Consider $n=[\frac{p}{q}]$ and divide the interval $[n,n+1]$ into $q$ 
sub-intervals $[\frac iq, \frac{i+1}{q}]$ for $i= nq, \ldots ,(n+1)q$. The 
$q$ numbers $\frac{j}{q-1}$ for $j=n(q-1), \ldots ,(n+1)(q-1)$ must divide 
among the $q$ sub-intervals, and there can be at most one in each 
interval because $\frac{1}{q} < \frac{1}{q-1}$.
\end{proof}

\begin{dfn}
We will call $(M \setminus (U_1 \cup U_2 \cup U_3), \ \xi|_{M \setminus (U_1 \cup U_2 \cup U_3)})$ the 
{\em background} of $(M,\ \xi)$.
\end{dfn}

The possible backgrounds will be studied in the next section.

\begin{lemma}\label{bound1}
Let $\xi$ be a tight contact structure with $t(\xi)<0$ on a Seifert manifold
$M(e_0; r_1, r_2, r_3)$ with integer Euler class $e_0=-1$. Then $t(\xi)<-1$.
\end{lemma}
\begin{proof}
Assume by contradiction that $t(\xi)=-1$.
By Proposition \ref{decomposizione} $\frac{p_i}{q}=-1$.
We would like to apply the classification theorem 
\cite[Lemma 5.1]{honda:2} to the background of $(M, \xi)$. However Honda 
orients the boundary by the outward normal convention, and uses 
the bases coming from the product structure, therefore with his 
conventions the boundary slopes of the background become $1$, $1$, and 
$0$.
By \cite[Lemma 5.1]{honda:2} point 3(a) the background of $(M, \xi)$ has
a vertical Legendrian curve $L$ with $tb(L)=0$, 
contradicting the hypothesis $t(\xi)<0$.  
\end{proof}

The following lemma is a technical observation which will be 
repeatedly useful in the paper.
\begin{lemma} \label{exploitmaximality}
Let $L$ be a maximally twisting Legendrian curve. If $A$ is a convex 
annulus with Legendrian boundary and one of its boundary components 
coincides with $L$, then the dividing set of $A$ contains no arcs
with both endpoints on $L$.
\end{lemma}
\begin{proof}
A dividing curve on $A$ with both endpoints on $L$ gives a bypass 
attached to $L$ as explained in \cite[Proposition 3.17]{honda:1}.
It is well known that the attachment of a bypass decreases the twisting
number (see 
\cite[Lemma 2.20]{etnyre:0})  contradicting our 
assumption.
\end{proof}

\begin{prop} \label{pendenze}
Let $M=M(e_0; r_1, r_2, r_3)$ be a small Seifert manifold with integer 
Euler class $e_0=-1$ or $e_0=-2$, and let $\xi$ be a tight contact 
structure on $M$ with maximal twisting $t(\xi)=-q<0$. If $\frac{p_1}{q}$, 
$\frac{p_2}{q}$, and $\frac{p_3}{q}$ are  the boundary slopes of the 
background of $(M, \xi)$, then $p_1+p_2+p_3=e_0q-1$. 
\end{prop}
\begin{proof}
If $t(\xi)=-1$, then $e_0=-2$ by Lemma \ref{bound1}. In this case 
Proposition \ref{decomposizione} implies $p_1=p_2=p_3=-1$, so the equality 
holds. 

Assume now $t(\xi)=-q < -1$.
Let $A$ be a convex vertical annulus between $- \partial (M \setminus U_1)$ and 
$- \partial (M \setminus U_2)$. The dividing set of $A$ has no boundary parallel 
dividing 
curves by Lemma \ref{exploitmaximality} because $-q$ is the maximal 
twisting number.  

By the edge-rounding lemma \cite[Lemma 3.11]{honda:1}, a neighbourhood
of $U_1 \cup U_2 \cup A$  has boundary slope $\frac{p_1+p_2+1}{q}$. A complement 
of this neighbourhood is a tubular neighbourhood $U_3'$ of $F_3$ 
containing $U_3$ such that $- \partial (M \setminus U_3')$ has slope $-\frac{p_1+p_2+1}{q} 
+e_0= -\frac{p_1+p_2+1-q e_0}{q}$. 
If $-\frac{p_1+p_2+1-q e_0}{q} \neq \frac{p_3}{q}$, then by  
Lemma \ref{algebrico} there would be a rational number 
$\frac{p'}{q-1} \in [\frac{p_3}{q}, -\frac{p_1+p_2+1}{q}]$. Then by 
\cite[Proposition 4.16]{honda:1}
 there would be a convex torus parallel to $- \partial (M \setminus U_3')$ in $U_3' \setminus U_3$ 
with slope $\frac{p'}{q-1}$ computed with respect to the basis of 
$- \partial (M \setminus U_3)$. This would contradicts $t(\xi)=-q$.
\end{proof}

\begin{cor} \label{gap}
Let $\xi$ be a tight contact structure with $t(\xi)<-1$ on a small Seifert 
manifold $M=M(e_0; r_1, r_2, r_3)$ with integer Euler class $e_0=-2$. Then 
$t(\xi) \leq -4$.
\end{cor}
\begin{proof}
 Let $\frac{p_1}{q}$, $\frac{p_2}{q}$, and $\frac{p_3}{q}$ be the boundary 
slopes of the background of $(M, \xi)$ as in Proposition 
\ref{decomposizione}. 
Recall that Proposition \ref{decomposizione}(2) implies 
$-1 < \frac{p_i}{q}$ (which is equivalent to $-q-p_i \leq -1$) when $q>1$.
From Proposition \ref{pendenze} we get $p_1=-p_2-p_3-2q-1$, therefore 
$p_1 \leq -3$. This implies $q>3$ because $\frac{p_1}{q} > -1$.  
\end{proof}

\section{Tight contact structures on $\Sigma \times S^1$.}
\begin{figure}\centering
\includegraphics[width=4cm]{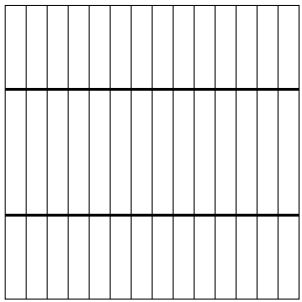}
\caption{The standard characteristic foliation ${\mathcal F}(0,1)$ on 
$T^2$.}
\label{standard.fig}
\end{figure}

In this section we classify all possible backgrounds. Given a pair of 
integer numbers $(p,q)$ with $q>0$, we denote by 
${\mathcal F}(p,q)$ the standard characteristic foliation on $T^2$ with 
vertical Legendrian ruling, slope $\frac pq$ and $2|(p,q)|$ Legendrian 
divides. Figure \ref{standard.fig} shows ${\mathcal F}(0,1)$.

\begin{lemma}\label{vertical}
Let $\Sigma$ be a pair of pants, and consider coordinates on 
$- \partial (\Sigma \times S^1)= T_1 \cup T_2 \cup T_3$ coming from the product structure.
For every triple of integer numbers $(p_1, p_2, q)$ with 
$q>0$ there exists a contact structure $\beta(p_1, p_2, q)$ on $\Sigma \times S^1$ 
which is tangent to the fibres, and induces the characteristic 
foliations ${\mathcal F}(p_1,q)$, ${\mathcal F}(p_2,q)$, and 
${\mathcal F}(p_3,q)$ with $p_3= -(p_1+p_2+1)$ on the three components of 
the boundary of $\Sigma \times S^1$.
\end{lemma}

\begin{figure} \centering
\psfrag{T1}{$T_1$}
\psfrag{T2}{$T_2$}
\psfrag{A}{$A \times S^1$}
\psfrag{B}{$B$}
\includegraphics[width=11cm]{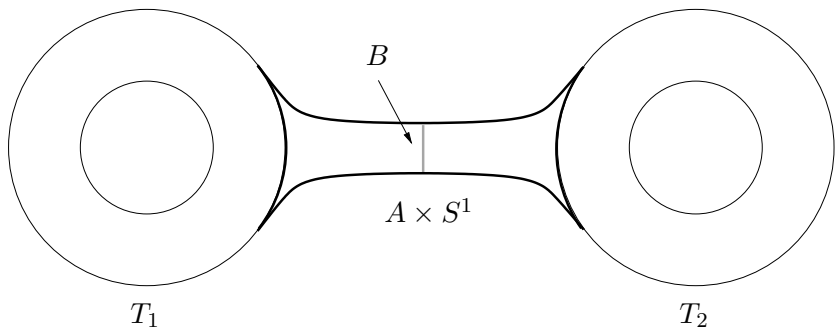}
\caption{The construction of $\beta(p_1,p_2,q)$.}
\label{raccordo.fig}
\end{figure}

\begin{proof}
Let $T_1 \times [0,1]$ and $T_2 \times [0,1]$ be invariant neighbourhoods of 
standard tori with characteristic foliations ${\mathcal F}(p_1,q)$ and
 ${\mathcal F}(p_2,q)$ respectively. Inside $\R^2 \times S^1$ with
the contact structure defined by the $1$--form $\cos (2 \pi qz)dx+
\sin (2 \pi qz)dy$ we take the subset $A \times S^1$, where $A$ is the subset 
of the plane portrayed in Figure \ref{raccordo.fig}. If we glue 
$T_1 \times [0,1]$, $T_2 \times [0,1]$, and $A \times S^1$ as in the picture, we obtain
$\Sigma \times S^1$ with the contact structure $\beta(p_1, p_2, q)$.
\end{proof}

\begin{lemma}
$\beta(p_1, p_2, q)$ is tight and its maximal twisting 
number is $-q$.
\end{lemma}
\begin{proof}
Suppose there is an overtwisted disc or a vertical Legendrian curve 
with twisting number greater than $-q$. Then there is a convex annulus 
$B'$ with Legendrian boundary which is smoothly isotopic to the 
annulus $B$ shown in Figure \ref{raccordo.fig} and is disjoint from 
the overtwisted disc or the Legendrian curve. By the Isotopy 
Discretisation Lemma \cite[Lemma 3.10]{gluing} 
%
% footnote
\footnote{The proof in the reference is incomplete because it assumes 
that non convex surfaces form a discrete set in any generic contact 
film. This is not true, however any contact film starting and ending 
with a convex surface can be isotoped relative to the boundary to a 
contact film with such property by \cite[Lemma 15]{giroux:6}.
}
there is a sequence of 
smoothly isotopic annuli $B=B_0, \ldots ,B_n=B'$, all with the same boundary, 
such that for each $i$ the interior parts of $B_{i-1}$ and $B_i$ are 
disjoint, and $B_i$ is obtained from $B_{i-1}$ by the attachment of a 
single bypass.

We will prove that all annuli $B_i$, and in particular $B_n=B'$, satisfy 
the following properties after a $C^0$--small modification:
\begin{enumerate}
\item $B_i$ is contact isotopic to $B_0$, and
\item if we cut $\Sigma \times S^1$ along $B_i$ and 
round the edges we obtain two connected components isomorphic to 
$T_1 \times [0,1]$ and $T_2 \times [0,1]$ with the invariant contact structures 
inducing the characteristic foliation ${\mathcal F}(p_1,q)$ on the 
boundary components of $T_1 \times [0,1]$, and ${\mathcal F}(p_2,q)$ on the 
boundary components of $T_2 \times [0,1]$.

\end{enumerate}
 This proves the Lemma because in the invariant contact structures on
$T_1 \times [0,1]$ and $T_2 \times [0,1]$ there are neither overtwisted discs, nor 
 vertical Legendrian curves with twisting number greater than $-q$.

Properties (1) and (2) hold for $B=B_0$ by construction. In order to 
prove that they hold for $B_i$, we assume that they hold for $B_{i-1}$, and 
prove that $B_i$ is contact isotopic to $B_{i-1}$. Since $B_i$ is disjoint 
from $B_{i-1}$ outside the boundary, it is contained in one of the two 
connected components of $\Sigma \times S^1 \setminus B_{i-1}$ (say in the one isomorphic to 
$T_1 \times [0,1]$ to fix the notation). The bypass sending $B_{i-1}$ to $B_i$ can
be trivial, can change the slope of $T_1 \times \{ 1 \}$ if $(p_1,q)=1$, or can
decrease the number of dividing curves of $T_1 \times \{ 1 \}$ if $(p_1,q)>1$. 
The last two options are impossible because the 
contact structure on $T_1 \times [0,1]$ is invariant. If the bypass sending 
$B_{i-1}$ to $B_i$ is trivial, then the dividing set of $B_i$ is isotopic to 
the dividing set of $B_{i-1}$, therefore we can modify $B_i$ in a small $C^0$
neighbourhood so that its 
characteristic foliation becomes isotopic to the characteristic 
foliation of $B_{i-1}$. After this modification $B_{i-1}$ is contact isotopic 
to $B_i$ and $\Sigma \times S^1 \setminus B_i$ is contact isotopic to $\Sigma \times S^1 \setminus B_{i-1}$ by the 
triviality of the attachment of trivial bypasses 
\cite[Lemma 2.10]{gluing}. The modification of the characteristic 
foliation  of $B_i$ is implicitly required by \cite[Lemma 2.10]{gluing}.
\end{proof}

\begin{prop}\label{unico}
Up to an isotopy not necessarily fixed on the boundary, 
$\beta(p_1,p_2,q)$ is the unique tight contact structure
on $\Sigma \times S^1$ with maximal twisting number $-q$ inducing the 
characteristic foliations ${\mathcal F}(p_1,q)$, ${\mathcal F}(p_2,q)$ and
${\mathcal F}(p_3,q)$ with $p_3=-(p_1+p_2+1)$ on the boundary.
\end{prop}

\begin{proof}
Take a convex annulus $A$ in $(\Sigma \times S^1, \xi)$ with Legendrian boundary 
between two vertical Legendrian ruling curves of $T_1$ and $T_2$ with 
twisting number $-q$. The dividing set of $A$ consists of curves with 
endpoints on different boundary components because of Lemma 
\ref{exploitmaximality} and the maximality of the twisting number of 
$\partial A$. After an isotopy of $\beta(p_1,p_2,q)$ not fixed on 
the boundary, we can assume that the dividing set of $A$ consists of 
$2q$ horizontal curves. After cutting along $A$ and rounding the edges 
by \cite[Edge Rounding]{honda:1} we obtain a toric annulus $N$
diffeomorphic to $T_3 \times [0,1]$ with characteristic foliation ${\mathcal 
F}(p_3,q)$ on both boundary components.
 The contact structure $\beta(p_1,p_2,q))$ restricted to $N$ is non-rotative, 
otherwise there would be a curve isotopic to a fibre with 
twisting number $0$, contradicting 
$t(\beta(p_1,p_2,q))=-q$. For the same reason there can be
no intermediate convex torus with fewer dividing curves.

Put the characteristic foliation on the boundary of $N$ in standard 
form so that each Legendrian ruling curve intersects each dividing 
curve exactly once. The considerations above imply that a Legendrian 
ruling curve of $\partial N$ minimises the twisting number in its isotopy 
class in $N$. Let $B$ be a convex annulus between two Legendrian 
ruling curves in different boundary components of $N$. By 
\cite[Theorem 2.2(4)]{honda:1} its dividing set determines 
the isotopy class of the contact structure on $N$. The dividing 
set of $B$ contains no boundary parallel dividing arcs because of 
Lemma \ref{exploitmaximality}, therefore it can be made horizontal 
with an isotopy of $N$ not fixed on the boundary, therefore 
  the contact structure on $N$ is isotopic to 
the invariant contact structure by the same isotopy.  This implies that
 $\beta(p_1,p_2,q)$ is unique 
up to isotopy not fixed on the boundary because it is determined by 
its restriction to the complement of a neighbourhood of $A$, and by 
the dividing set of $A$. 
\end{proof}

Since a contact isotopy of the background, even if it is not constant 
on the boundary, can be extended to a contact
isotopy of $(M, \xi)$ --- see \cite[Lemma 4.4]{ghiggini:1} --- we 
have the following corollary.
\begin{cor}
$(M, \xi)$ has a background which is isotopic to $\beta(p_1,p_2,q)$ with 
$q=-t(\xi)$.
\end{cor}

\section{Transverse contact structures}
In the following a {\em transverse contact structure} on a Seifert
 manifold will be a contact structure which is positively 
transverse to the Seifert fibration. This condition has strong 
consequences both on the contact structure and on the topology of the
underlying manifold. For example:

%Results analogous to those in this section have been proved 
%originally by Lisca and Mati\'c using different techniques. In 
%particular, they have 
%numerical conditions on the Seifert coefficient which determine if a 
%given small Seifert manifold admits a transverse contact structure. 
%The conditions we prove in Theorem \ref{costruzione} is different, 
%but of course equivalent, to Lisca and Mati\'c's one.

\begin{theorem} \label{tight}
(\cite[Theorem 2.2]{lisca-matic:transverse})
Transverse contact structures on Seifert manifolds are universally 
tight.
\end{theorem}

\begin{prop}\label{pippo}
If $\xi$ is a transverse contact structure on a small Seifert manifold 
$M$, then $t(\xi)<0$.
\end{prop}

\begin{proof}
The base space of a small Seifert manifold is a $2$--dimensional 
orbifold
with three cone points having $S^2$ as underlying surface. The three 
cone points are the images of the singular fibres. Any such orbifold 
is finitely covered in the sense of orbifolds by a smooth surfaces 
$\Sigma'$ (see \cite{scott}). 
The pull back of the Seifert fibration to $\Sigma'$ is a honest circle 
bundle $M'$, because the singular set of $\Sigma'$ is empty. If $\xi$ is 
transverse to the Seifert 
fibration of $M$, the pulled back contact structure $\xi'$ is transverse 
to the fibres of $M'$. By \cite[Theorem 2.3]{giroux:3}, a universally 
tight contact structure on a circle bundle over a 
surface is transverse if and only if there is no Legendrian curve with 
twisting number $0$ isotopic to a fibre.
\end{proof}

 \begin{rem}
A second proof of this lemma can be given by applying the slice 
Thurston--Bennequin inequality of Kronheimer and Mrowka to the 
regular fibres in the symplectic fillings constructed by Lisca and 
Mati\'c in the proof of \cite[Theorem 2.2]{lisca-matic:transverse}.
\end{rem}

\begin{lemma}\label{horizontal}
For every triple of integer numbers $(p_1, p_2, q)$ with $q>0$ we can 
perturb $\beta(p_1, p_2, q)$ in any arbitrarily small 
$C^{\infty}$--neighbourhood, and obtain a new contact 
structure $\tilde{\beta}(p_1, p_2, q)$ on $\Sigma \times S^1$ which is
transverse to the $S^1$--fibres and has convex boundary with the same 
dividing set as $\beta(p_1, p_2, q)$.
\end{lemma}
\begin{proof}
Fix a contact form $\alpha$ for $\beta(p_1, p_2, q)$. Let $dz$ 
be the pull-back of a volume form on $S^1$ to $\Sigma \times S^1$, and 
$\frac{\partial}{\partial z}$ be a vector field tangent to the $S^1$--fibration such 
that $dz(\frac{\partial}{\partial z}) =1$. Then, for $\epsilon$ small, the $1$--form
 $\alpha + \epsilon dz$ defines a 
contact structure which is transverse to the fibration because
$\alpha(\frac{\partial}{\partial z})=0$.

Put coordinates $(x,y,z)$ near a boundary component of 
$\Sigma \times S^1$, so that $z$ is the direction of the fibres and 
$x$ is the direction of the inward normal. The contact structure in a 
neighbourhood of a boundary component with characteristic foliation
${\mathcal F}(p,q)$ --- the type of characteristic foliation induced 
by $\beta$ on the boundary of $\Sigma \times S^1$ --- is locally defined by the 
$1$--form 
$\alpha= \cos(2 \pi (qz+py))dx + \sin (2 \pi (qz+py))dy$, and $\frac{\partial}{\partial x}$ is 
a contact vector field. It is straightforward to check that 
$\frac{\partial}{\partial x}$ is a contact vector field for the kernel of $\alpha + \epsilon dz$ 
too, and that the dividing set remains unchanged.
\end{proof}

\begin{lemma}\label{connection}
Let ${\mathcal F}_0$ and ${\mathcal F_1}$ be foliations on $T^2$ divided 
by the same multicurve $\Gamma$, and which are both transverse to the same 
choice of a vertical direction on $T^2$. Then there exists a transverse 
contact structure on $T^2 \times I$ which induces ${\mathcal F}_0$ on 
$T^2 \times \{ 0 \}$ and ${\mathcal F}_1$ on $T^2 \times \{ 1 \}$ as characteristic 
foliations.
\end{lemma}
\begin{proof}
Let $\omega$ be an area form on $T^2$. In the proof of 
\cite[Proposition II.3.6]{giroux:1}, Giroux constructs 
a family of vector fields $Y_s=(1-s)Y_0+sY_1$ and functions $v_s \colon T^2 \to \R$ 
such that $Y_i$ directs ${\mathcal F}_i$ for $i=0,1$, and $\iota_{Y_s}\omega +v_s dt$ 
is a contact form on $T^2 \times \R$ for all $s \in [0,1]$. By 
\cite[Lemma 2.3]{giroux:4}
there is a function $u \colon T^2 \times \R \to \R$ such that 
$\iota_{Y_t}\omega +u_t dt$ is a contact form on $T^2 \times [0,1]$. The transversality 
condition on ${\mathcal F}_0$ and ${\mathcal F}_1$ is equivalent to the 
existence of a $1$--form $\lambda$ on $T^2$ such that $\lambda(Y_i)>0$ for $i=0,1$. 
Since $\lambda(Y_t)>0$ for any $t \in [0,1]$ by linearity, the contact structure 
defined by $\iota_{Y_t}\omega +u_t dt$ is transverse.
\end{proof}

\begin{theorem} \label{costruzione}
A small Seifert manifold $M=M(e_0; r_1, r_2, r_3)$ with $e_0= -2, -1$ admits 
a positively
transverse contact structure if and only if there are integer 
numbers $p_1$, $p_2$, $p_3$ and $q>0$ such that
\begin{enumerate} 
\item $\frac{p_i}{q} < -r_i$, and
\item $p_1 +p_2 +p_3 =qe_0 -1$.
\end{enumerate}
\end{theorem}

\begin{rem}
Different, but equivalent, necessary and sufficient conditions for 
the existence of transverse contact structures on small Seifert 
manifolds have been proved by Lisca and Mati\'c using $4$--dimensional 
techniques and pre-existing results on taut foliations on small 
Seifert manifolds; see \cite{lisca-matic:transverse}.
\end{rem}

\begin{proof}
The ``only if'' part follows from Proposition \ref{decomposizione}, 
Proposition \ref{pendenze}, and Proposition \ref{pippo}.
For the ``if'' part take neighbourhoods $U_1$, $U_2$, and $U_3$ of the 
singular fibres as in Proposition \ref{decomposizione}, and for 
$i=1,2,3$ denote by $s_i$ the boundary slope of $U_i$ corresponding to 
$\frac{p_i}{q}$ in the basis of $- \partial (M \setminus U_i)$. Let $\xi_i$ be the contact 
structure on $U_i \cong D^2 \times S^1$ with coordinates $(\rho, \phi, \theta)$ defined 
by the equation $\xi_i = \ker (\cos (k_i \rho) d \theta + k_i \rho \sin (k_i \rho) d \phi)$, where 
$k_i$ has been chosen so that $\partial U_i$ is pre-Lagrangian and has slope $s_i$.
For $i=1,2,3$, $\xi_i$ is transverse to the Seifert fibration of $M$
restricted to $U_i$. 
In fact, by the invariance of $\xi_i$ in the directions $\theta$ and $\phi$, a  
tangency between $\xi_i$ and a regular fibre would produce a 
pre-Lagrangian torus where the leaves of the characteristic foliation 
coincides with the fibres of the Seifert fibration. A convex 
perturbation of this torus (as constructed in 
\cite[Lemma 3.4]{ghiggini:1}) would have infinite slope in the basis
of $- \partial (M \setminus U_i)$. 
This is impossible because the slopes in $U_i$, computed in the basis of 
$- \partial (M \setminus U_i)$, belong to the interval $[\frac{p_i}{q}, -r]$.

The boundary of $U_i$ can be made convex with $(p_i, q)$ dividing 
curves by a $C^{\infty}$-small perturbation of $\xi_i$, which therefore does not 
affect transversality; see \cite[Lemma 3.4]{ghiggini:1}. Then for 
$i=1,2,3$
we use Lemma \ref{connection} to make the characteristic foliations of 
$\tilde{\beta}(p_1,p_2,q)$ on $- \partial (M \setminus U_i)$ and the characteristic foliation 
of $\xi_i$ on $\partial U_i$ match under the gluing maps, still without affecting 
transversality. When we glue all pieces together, we get
a contact structure on $M$ which is always transverse to the Seifert 
fibration.
\end{proof}

\begin{rem}
If $e_0=-2$ we can always take $p_1=p_2=p_3=-1$ and $q=1$, so every small 
Seifert manifold with $e_0=-2$ admits a transverse tight contact 
structure.
\end{rem}

\begin{cor} 
A small Seifert manifold $M(e_0; r_1,r_2,r_3)$ 
admits a tight contact structure $\xi$ with $t(\xi)<0$ if and only if it 
admits a transverse contact structure.
\end{cor}
\begin{proof}
When $e_0 \neq -1$ it follows from works of Wu \cite{hwu} and Lisca and 
Mati\'c \cite{lisca-matic:transverse} combined.
 When $e_0 \neq -1$
one direction has been proved in Proposition \ref{pippo}. To prove the 
other direction assume that $M(e_0; r_1, r_2, r_3)$ admits a tight contact 
structure with negative twisting. Then, combining Proposition 
\ref{decomposizione}(1) with Proposition \ref{pendenze}, we obtain
integer number $p_1,p_2,p_3$, and $q>0$ such that $(p_1,q)=1$, 
$\frac{p_i}{q}< -r_i$, and $p_1 + p_2 + p_3 =qe_0 -1$, therefore 
$M(e_0; r_1, r_2, r_3)$ admits a transverse contact structure by Theorem 
\ref{costruzione}.
\end{proof}

\begin{proof}[Proof of Theorem \ref{principale}]
Theorem \ref{principale} follows from Proposition \ref{decomposizione},
Proposition \ref{pippo}, and Theorem \ref{costruzione}.
\end{proof}

\section{Tight contact structures on $L$--spaces}

In this section we classify tight contact structures on those small 
Seifert manifolds with $e_0=-2$ which are $L$--spaces. The condition of
being an $L$--space is used to give a bound of the maximal twisted 
number, which in turn gives  
 an upper bound 
on the number of tight contact structures. Finally we construct 
enough distinct tight contact structures to match the upper bound. 
After the bound on the maximal twisted number (Proposition \ref{L}) 
the proof goes on like in \cite{hwu}.

\begin{prop} \label{L}
Let $M=M(e_0; r_1, r_2, r_3)$ be a small Seifert manifold with integer 
Euler class $e_0= -2$. If $M$ admits a tight contact structure $\xi$ with 
$t(\xi)< -1$, then $-M$ admits a positively transverse contact structure.
 \end{prop}

\begin{proof}
Assume that $M$ admits a tight contact structure $\xi$ with $t(\xi)<-1$.
Then by corollary \ref{gap} $t(\xi)<-4$.
Take numbers $p_1$, $p_2$, $p_3$, and $q$ as in Proposition 
\ref{decomposizione} with $q>4$ because $t(\xi)<-4$.
For $i=1,2,3$ we have 
$\frac{p_i}{q} < -r_i < \frac{p_i+2}{q-2}$ because $\frac{p_i}{q}$ is 
negative. Also, $p_1+p_2+p_3=-2q-1$.

Define $p_i'= -(p_i+q)$ and $q'=q-2$ so that $\frac{p_i'}{q'}= -1- 
\frac{p+2}{q-2}$, then
$\frac{p_i'}{q'} < -1+r_i$ and $p_1'+p_2'+p_3'=-q'-1$. This implies that $-M$
carries a transverse contact structure because $-M(-2; r_1,r_2,r_3)=
M(-1, 1-r_1, 1-r_2, 1-r_3)$.

\end{proof}
 
By using \cite[Theorem 1.1]{lisca-stipsicz:6} we have the following 
corollary
\begin{cor}\label{finalmente}
If a small Seifert manifold $M$ with integer Euler class $e_0 \leq -2$ is
an $L$--space, then $t(\xi)= -1$ for any tight contact structure $\xi$ on
$M$.
\end{cor}

\begin{proof}[Proof of Theorem \ref{classificazione}]
The part of the statement concerning the maximal twisted number 
%has been proved in Proposition \ref{L}. 
is Corollary \ref{finalmente}. It implies, by 
Proposition \ref{decomposizione} and Proposition \ref{unico}, that the
only possibility for the background of $\xi$ is $\beta(-1,-1,1)$. After some 
easy arithmetics, from the classification of tight contact structures 
on solid tori \cite[Theorem 2.3]{honda:1} we have $T(r_1)$, $T(r_2)$, and 
$T(r_3)$ possible isotopy classes of tight contact structures on the 
neighbourhoods of the singular fibres $U_1$, $U_2$, and $U_3$ respectively.
 Altogether, they give an upper bound of at most $T(r_1)T(r_2)T(r_3)$ 
isotopy classes of tight contact structures on $M(-2; r_1, r_2, r_3)$.

In order to construct $T(r_1)T(r_2)T(r_3)$ non isotopic tight contact 
structures we perform Legendrian surgery on all possible Legendrian 
realisations of the link in Figure \ref{link.fig} with the appropriate 
Thurston--Bennequin numbers of the components. We have a unique 
possibility for the central unknot with surgery coefficient $-2$: it 
must be a Legendrian unknot with Thurston--Bennequin number $-1$. An 
unknot in one of the ``legs'' with surgery coefficient $a^{(i)}_j$ must be 
made Legendrian with Thurston--Bennequin number $a^{(i)}_j+1$, therefore we 
have $|a^{(i)}_j+1|$ choices for its rotation number. Varying over all 
possible choices of the rotation number for all components of the 
links produces $T(r_1)T(r_2)T(r_3)$ Stein fillable contact structures 
whose fillings are all diffeomorphic, but have pairwise distinct first 
Chern classes. By \cite[Theorem 2]{plam:1} all the contact 
structures constructed from these surgeries have pairwise distinct and 
non trivial Ozsv\'ath-Szab\'o invariants with coefficients in $\Z / 2\Z$, 
therefore they are pairwise non isotopic (see also 
\cite{lisca-matic:2} for a similar way to prove the same result). 
Since $M$ is an $L$--space by \cite[Theorem 1.1]{lisca-stipsicz:6}, 
for any ${\rm Spin}^c$--structure $\mathfrak{s}$ there is only one non 
zero element in $\widehat{HF}(-M, \mathfrak{s})$ with $\Z / 2\Z$ 
coefficients. This implies that the tight contact structures defined 
by all possible Legendrian surgeries on the link in Figure 
\ref{link.fig} define pairwise distinct ${\rm Spin}^c$--structures.  
\end{proof}

\bibliographystyle{plain}
\bibliography{contatto}
\end{document}